\documentclass{article}
\usepackage{amsmath}
\usepackage{amssymb}
\usepackage{a4wide}

\newcommand\bX{\mathbf X}

\newcommand\X{\mathbf X}
\newcommand\bR{{\mathbb R}}

\newcommand\bN{{\mathbb N}}

\newcommand\bI{I}
\newcommand\e{{\mathbb E}}
\newcommand\p{{\mathbb P}}

\newcommand\tr{{\mathrm{tr}}}
\newcommand\cF{{\mathcal F}}

\newtheorem{theorem}{Theorem}[section]%
\newtheorem{corollary}[theorem]{Corollary}%
\newtheorem{lemma}[theorem]{Lemma}%
\newtheorem{proposition}[theorem]{Proposition}%
\newtheorem{remark}[theorem]{Remark}%

\begin{document}

\begin{center}
\Large Lower bounds on the smallest eigenvalue of a sample covariance matrix.  
\end{center}
\begin{center}
\large Pavel~Yaskov\footnote{Steklov Mathematical Institute of RAS, Russia\\
 e-mail:yaskov@mi.ras.ru\\Supported
    by RNF grant 14-21-00162 from the Russian Scientific Fund.}
 \end{center}

\begin{abstract}We provide tight lower bounds on the smallest eigenvalue of a sample covariance matrix of a centred isotropic random vector under weak or no assumptions on its components.
\end{abstract}

\begin{center}
{\bf Keywords:} Covariance matrices; Gram matrices; Random matrices. 
\end{center}
\section{Introduction} Lower bounds on the smallest eigenvalue of a sample covariance matrix (or a Gram matrix) play a crucial role in the least squares problems in high-dimensional statistics (see, for example, \cite{O}). These problems motivate the present work.

For a random vector $X_p$ in $\mathbb R^p$,
consider a random $p\times n$ matrix $\X_{pn}$ with  independent columns $\{X_{pk}\}_{k=1}^n$ distributed as $X_p$ and the Gram matrix
 \[\X_{pn}\X_{pn}^\top=\sum_{k=1}^n X_{pk}X_{pk}^\top.\]
 If $X_p$ is centred, then $n^{-1}\X_{pn}\X_{pn}^\top$ is the sample covariance matrix corresponding to the random sample $\{X_{pk}\}_{k=1}^n$. For simplicity, we will further assume that $X_p$ is isotropic, i.e. $\e X_pX_p^\top=I_p$ for a $p\times p$ identity matrix $I_p$, and consider only those  $p$ which are not greater than $n$  (otherwise $\X_{pn}\X_{pn}^\top$ would be degenerate).  
 
  In this paper we derive sharp lower bounds for $\lambda_p(n^{-1}\X_{pn}\X_{pn}^\top)$, where $\lambda_p(A)$ is the smallest  eigenvalue of a  $p\times p$ matrix $A$. We try to impose as few restrictions on the components of $X_p$ as possible. In proofs we use the same strategy as in \cite{SV}.

\section{Main results}
Put $c_p(a)=\inf \e \min\{(X_p,v)^2,a\}$, $C_p(a)=\sup\e (X_p,v)^2 \min\{(X_p,v)^2,a\}$, 
\[L_p(\alpha)=\sup\e |(X_p,v)|^{2+\alpha}
\quad\text{and}\quad K_p=\inf \e|(X_p,v)|\]
for given $a,\alpha>0,$ where all suprema and infima are taken over $v\in\bR^p$ with $\|v\|=1,$ and $\|v\|=(\sum_{i=1}^pv_i^2)^{1/2}$ is the Euclidean norm of $v=(v_1,\ldots,v_p).$  Denote also by $M_p(\alpha)$ the infimum over all $M>0$ such that 
\[\p(|(X_p,v)|>t)\leqslant\frac{M}{t^{2+\alpha}}\quad\text{for all $t>0$ and  $v\in\bR^p,$ $\|v\|=1.$}\]
Our main lower bounds are as follows.
\begin{theorem}\label{thm1}
If $X_p$ is an isotropic random vector in $\mathbb R^p$ and $p/n\leqslant y$ for  some $y\in(0,1)$,  then, for all 
$a>0$,
\[\lambda_p(n^{-1}\X_{pn}\X_{pn}^\top)\geqslant c_p(a)-\frac{C_p(a)}{a} -5a y+\frac{\sqrt{C_p(2a)}Z}{\sqrt{n}}\]
for a centred random variable $Z=Z(p,n,a)$ with  $\p(Z<-t)\leqslant e^{-t^2/2}$, $t>0$.
\end{theorem}

\begin{theorem}\label{thm2}
Let $X_p$ be an isotropic random vector in $\mathbb R^p$, $p/n\leqslant y$ for  some $y\in(0,1)$. If $L_p(2)<\infty$, then 
\[\lambda_p(n^{-1}\X_{pn}\X_{pn}^\top)\geqslant 1-4C\sqrt{y}+\frac{CZ}{\sqrt{n}}\]
for $C= \sqrt{L_p(2)}$ and some $Z=Z(p,n)$   with $\e Z=0$ and $\p(Z<-t)\leqslant e^{-t^2/2}$, $t>0$. Moreover, there are universal constants $C_0,C_1,C_2>0$ such that
\[\lambda_p(n^{-1}\X_{pn}\X_{pn}^\top)\geqslant C_0K_p^2+\frac{C_1Z}{\sqrt{n}}\] whenever
$y\leqslant C_2K_p^2$ and $Z=Z(p,n)$ as above.
\end{theorem}

Useful bounds for $c_p(a)$ and $C_p(a)$ in terms of $L_p(\alpha)$ and $M_p(\alpha)$ are given in the following proposition.

\begin{proposition}\label{p1} 
Let $X_p$ be an isotropic random vector in $\mathbb R^p.$ Then, for all $a,\alpha>0,$
\[c_p(a)\geqslant 1-\frac{L_p(\alpha)}{a^{\alpha/2}}\quad\text{and}\quad c_p(a)\geqslant1-\frac{2\alpha^{-1}M_p(\alpha)}{a^{\alpha/2}}.\]
 In addition, for all $\alpha\in(0,2]$ and each $a>0$, $C_p(a)$ is bounded from above
 by
\[ a^{1-\alpha/2} L_p(\alpha)\quad\text{and}\quad  (1+2/\alpha)M_p(\alpha)a^{1-\alpha/2}+\begin{cases}2M_p(\alpha)a^{1-\alpha/2}/(1-\alpha/2),&\alpha\in(0,2),\\
2M_p(2)\log\max\{a,1\}+1,&\alpha=2.\end{cases}\]  
\end{proposition}

\section{Applications}
We now describe different corollaries of Theorem \ref{thm1} and Theorem \ref{thm2}.  
The next corollary extends Theorem 1.3 in \cite{KM} and Theorem 3.1 in \cite{O} (for $A_i=X_{pi}X_{pi}^\top$).  
\begin{corollary}\label{cl1} 
Let $X_p$ be an isotropic random vector in $\mathbb R^p$, $p/n\leqslant y$ for some $y\in(0,1)$ and $L_p(\alpha)<\infty$ for some $\alpha\in(0,2]$. Then, with probability at least $1-e^{-p}$,
\[\lambda_p(n^{-1}\X_{pn}\X_{pn}^\top)\geqslant 1-C_\alpha y^{\alpha/(2+\alpha)},\]
where
\[C_\alpha=\begin{cases}9(L_p(\alpha))^{2/(2+\alpha)},&\alpha\in(0,2),\\
(4+\sqrt{2})\sqrt{L_p(2)},&\alpha=2.
\end{cases}\]
\end{corollary}\nopagebreak
\begin{remark} {\normalfont
One may  further weaken assumptions in Corollary \ref{cl1}. Namely, one may assume that  $M_p(\alpha)<\infty$ for some $\alpha\in(0,2).$ The conclusion of  Corollary \ref{cl1} will still hold with some $C_\alpha>0$ that depends only on $\alpha$ and $M_p(\alpha)$. In the case $\alpha=2$, one would have a lower bound of the form 
$1-C_2\sqrt{y \log(e/y)}$ with $C_2>0$ depending only on $M_p(2).$}
\end{remark}

Theorems \ref{thm1} and \ref{thm2}  improve Theorem 2.1 in \cite{SV} as the next corollary shows.

\begin{corollary}\label{cl2} 
Let $X_p$ be an isotropic random vector in $\mathbb R^p$.  If $L_p(\alpha)<\infty$ for some $\alpha\in(0,2)$ and  $p/n\leqslant \varepsilon^{1+2/\alpha}/(10(4L_p(\alpha))^{2/\alpha})$, then 
\[\e \lambda_p(n^{-1}\X_{pn}\X_{pn}^\top)\geqslant 1-\varepsilon.\]
The same conclusion holds if $L_p(2)<\infty$ and  $n\geqslant 16L_p(2)\varepsilon^{-2}p$.
\end{corollary}

Let us formulate the final corollary that improves Theorem 3.1 in \cite{KM} for small $K_p$.
\begin{corollary}\label{cl3} 
Let $X_p$ be an isotropic random vector in $\mathbb R^p$. Then there are universal constants $C_0^*,C_1^*,C_2^*>0$ such that, with probability at least  $1-\exp\{-C_1^*K_p^4 n\}$, 
\[\lambda_p(n^{-1}\X_{pn}\X_{pn}^\top)\geqslant C_0^* K_p^2\]
when $p/n\leqslant C_2^*K_p^2.$
\end{corollary}

The range of applicability of Corollary \ref{cl3} is very wide.  Namely, there exist some universal constant $K>0$  such that $K_p\geqslant K$ for a very large class of isotropic random vectors $X_p$. By Corollary \ref{cl3}, this means that   $\lambda_p(n^{-1}\X_{pn}\X_{pn}^\top)$ is separated from zero by an universal constant. 

The existence of $K$ follows from results related to Kashin's decomposition theorem. The infinite dimensional version of this theorem is given in Kashin \cite{K} (for a proof, see \cite{Kr}). It states the following.
\begin{quote}
There is an universal constant $K>0$ such that $L_2(0,1)=H_1\oplus H_2$ for some linear subspaces of $H_i\subset L_2(0,1),$ $i=1,2,$ such that $\|x\|_1\geqslant K\|x\|_2$  for all $x\in H_1\cup H_2,$ where $\|x\|_d$ is the standard norm in $L_d(0,1)$, $d=1,2$.
\end{quote}
Let $(\Omega,\cF,\p)$ be an underlying probability space. Assume that $\Omega=(0,1),$ $\cF$ is the Borel $\sigma$-algebra and $\p$ is the Lebesgue measure. If all components of $X_p=(x_1,\ldots,x_p)$ are in $ H_1$, or all components of $X_p$ are in $H_2$, then $K_p\geqslant K$. 

If we consider only discrete random vectors $X_p$, we may say more. Namely, Kashin \cite{K74} proved that, for any $\delta>0$ and all $N\in\bN$, $\bR^N$  contains a linear subspace $H$ with $\dim H\geqslant (1-\delta) N$ such that $|e|_1\geqslant K|e|_2$ for some $K=K(\delta)>0$ not depending on $N$ and all $e=(e_1,\ldots,e_N)\in H$,\footnote[1]{In fact, the Haar measure of such orthogonal matrices $C$ that $H=CH_1$ satisfies this property is greater than $1-2^{-N}$ for some $K=K(\delta)>0$, where $H_1=\{(e_1,\ldots,e_N)\in\bR^N:e_i=0,\;i\geqslant (1-\delta)N+1\}$ (see \cite{K74}).} where \[|e|_d=\Big(\frac{1}{N}\sum_{i=1}^N|e_i|^d\Big)^{1/d},\quad d=1,2.\]
In particular, if $\{e^{(k)}\}_{k=1}^p$ is {\it any} orthonormal system in $H$ and $\{x^{(i)}\}_{i=1}^N$ are columns of the $p\times N$ matrix with rows $\{(e^{(k)})^\top\}_{k=1}^p$,  then, for all $v=(v_1,\ldots,v_p)\in\bR^p$ with $\|v\|=\sqrt{\sum_{j=1}^pv_j^2}=1,$ 
\[K=K\Big(\frac{1}{N}\sum_{i=1}^N|(x^{(i)},v)|^2\Big)^{1/2}=K\Big|\sum_{k=1}^pv_k e^{(k)}\Big|_2\leqslant \Big|\sum_{k=1}^pv_k e^{(k)}\Big|_1=\frac1N\sum_{i=1}^N|(x^{(i)},v)|.\]
If $X_p$ is such that $\p(X_p=x^{(i)})=1/N,$ $1\leqslant i\leqslant N$, then $K_p\geqslant K=K(\delta).$

\section{Proofs.}

In proofs of Theorem \ref{thm1} and Theorem \ref{thm2}, we follow the strategy of Srivastava and Vershynin \cite{SV}. The key step is the following lemma.

 \begin{lemma}\label{l1}
 Let $A$  be a $p\times p$ symmetric matrix with $A\succ 0$, $v\in\bR^p$, $l\geqslant 0$, $\varphi>0$,
\begin{equation}
Q(l,v)=v^\top(A-l I_p)^{-1}v\quad\text{and}\quad
q(l,v)=\frac{v^\top(A-l I_p)^{-2}v}{\tr(A-l I_p)^{-2}},
\label{Qq}
\end{equation}
hereinafter $A\succ 0$ means that $A$ is positive definite. If $A-lI_p\succ 0,$ $\tr(A-lI_p)^{-1}\leqslant \varphi$ and
\[\Delta=\frac{q(l,v)}{1+3\varphi q(l,v)+Q(l,v)},\]
then $A-(l+\Delta )I_p\succ 0$  and  $\tr(A+vv^\top-(l+\Delta) I_p)^{-1}\leqslant \varphi$.
\end{lemma}
The proof of Lemma \ref{l1} is given in Appendix.  

The strategy itself consists in the following. Let $A_0$ be a $p\times p$ zero matrix and \[A_{k}=\sum_{j=1}^k X_{pj}X_{pj}^\top,\quad 1\leqslant k\leqslant n.\] Consider some $\varphi>0$ and  take $l_0=-p/\varphi$ that satisfies   $\tr(A_0-l_0I_p)^{-1}= \varphi$.

Put $l_k=l_{k-1}+\Delta_k$ for $1\leqslant k\leqslant n$,  where 
\[\Delta_k=\frac{q_k(l_{k-1},X_{pk})}{1+3\varphi q_k(l_{k-1},X_{pk})+Q_k(l_{k-1},X_{pk})},\] 
$Q_k(l_{k-1},X_{pk})$ and $q_k(l_{k-1},X_{pk})$ are defined as $Q(l ,v)$ and $q(l ,v)$ in \eqref{Qq} with $A=A_{k-1}$   and $v=X_{pk}$. Applying Lemma \ref{l1} iteratively, we infer that $\tr(A_k-l_kI_p)^{-1}\leqslant \varphi$ and $A_k-l_k I_p\succ 0$ for all $1\leqslant k\leqslant n$. Therefore, \[\lambda_p(\bX_{pn}\bX_{pn}^\top)=\lambda_p(A_n)\geqslant l_n=l_0+\Delta_{1}+\ldots+\Delta_n.\]

Let $\e_{k}=\e(\,\cdot\,|X_{p1},\ldots,X_{pk})$, $1\leqslant k\leqslant n$, and  $\e_0=\e$. We have \begin{equation}\label{star}
\lambda_p(n^{-1}\bX_{pn}\bX_{pn}^\top)\geqslant -\frac{p}{n\varphi}+\frac{1}{n} \sum_{k=1}^n\e_{k-1}\Delta_k+\frac{Y}{\sqrt{n} },\end{equation}
 where  $Y=n^{-1/2}\sum_{k=1}^n(\Delta_k-\e_{k-1}\Delta_k).$ 

To apply estimate \eqref{star}, we need to choose $\varphi$ and obtain good lower bounds for $\e_{k-1}\Delta_k$ as well as upper bounds for $\p(Y<-t),$ $t<0$. The next lemmata which proofs are given in Appendix provide such bounds.
 
\begin{lemma}\label{lm1}
Let $U$ and $V$ be non-negative random variables. Then, for all $a>0$, 
\[\e \frac{U}{1+V}\geqslant 
\frac{|\e\min\{U,a\}|^2}{\e\min\{U,a\}+\e V \min\{U,a\}}.\]
In addition, if $\e U=1,$ then  $\e U/(1+V)\geqslant 1/(1+\e UV)$. Moreover, 
\[\e\frac{U}{1+V}\geqslant \frac{|\e\sqrt{U}|^2}{1+\e V}. \]
\end{lemma}

\begin{lemma}\label{lm2}
Let $X_p$ be an isotropic random vector in $\bR^p$, $A,B\succ 0$ be a $p\times p$ symmetric matrices with $\tr(A)=1$ and $\tr(B)\leqslant 1$ that are  simultaneously diagonalisable.  If
\[\Delta=\frac{X_p^\top AX_p}{1+b^{-1}(X_p^\top AX_p+X_p^\top B X_p/3)}\]
for some $b>0$, then, for any $a>0$, \[\e\Delta\geqslant c_p(a)-\frac{5C_p(a)}{3b}\quad\text{and}\quad \e \Delta^2\leqslant C_p(b).\]
In addition, if $L_p(2)<\infty$, then $\e\Delta\geqslant 1-4L_p(2)b^{-1}/3$ and $\e \Delta^2\leqslant L_p(2)$. Moreover, \[\e \Delta\geqslant \frac{K_p^2}{1+4(3b)^{-1}}.\] 
\end{lemma}

\begin{lemma}\label{lm3}
Let $(D_k)_{k=1}^n$ be a sequence of non-negative random variables adapted to a filtration $(\cF_k)_{k=1}^n$ such that $\e(D_k^2|\cF_{k-1})\leqslant 1$ a.s. for $k=1,\ldots,n$, where $\cF_0$ is the trivial $\sigma$-algebra. If 
\[Z=\frac{1}{\sqrt{n}}\sum_{k=1}^n(D_k-\e(D_k|\cF_{k-1})),\]
then $\p(Z<-t)\leqslant \exp\{-t^2/2\}$ for all $t>0$.
\end{lemma}

{\bf Proof of Theorem \ref{thm1}.}
Take in Lemma \ref{lm2} $X_p=X_{pk},$
\begin{equation}\label{def}
A=\frac{(A_{k-1}-l_{k-1} I_p)^{-2}}{\tr(A_{k-1}-l_{k-1} I_p)^{-2}},\quad B=(A_{k-1}-l_{k-1} I_p)^{-1}/\varphi,\quad a=\frac1{5\varphi},\quad
b=\frac{5a}3=\frac{1}{3\varphi}.
\end{equation}
Clearly $A$ and $B$ commute hence they are simultaneously diagonalizable. Additionally, we have
$\tr(A)=1$ and $\tr(B)=\tr (A_{k-1}-l_{k-1}I_p)^{-1}/\varphi\leqslant1$. Using Lemma \ref{lm2}, we arrive at the lower bounds
\[\e_{k-1}\Delta_k\geqslant c_p(a)-\frac{C_p(a)}a,\quad 1\leqslant k\leqslant n,\]
hereinafter all inequalities with conditional mathematical expectations hold almost surely. By  \eqref{star}, the latter  implies that
\[\lambda_p(n^{-1}\bX_{pn}\bX_{pn}^\top)\geqslant  c_p(a)-\frac{C_p(a)}a-\frac{5ap}{n}+\frac{\sqrt{C_p(2a)}Z}{\sqrt{n}},\]
where  \[Z=\frac{1}{\sqrt{C_p(2a)n}}\sum_{k=1}^n(\Delta_k-\e_{k-1}\Delta_k).\]
 Note that $(\Delta_k-\e_{k-1}\Delta_k)_{k=1}^n$ is a martingale difference sequence  with respect to the natural filtration of $(X_{pk})_{k=1}^n$. Obviously, $\e Z=0$.  By Lemma \ref{lm2}, $\e_{k-1}\Delta_k^2\leqslant C_p(b)\leqslant C_p(2a)$. Therefore, Lemma \ref{lm3} with $D_k=\Delta_k/\sqrt{C_p(2a)}$ yields that $\p(Z<-t)\leqslant \exp\{-t^2/2\},$ $t>0.$
Thus we have proven Theorem \ref{thm1}.
\\\\
{\bf Proof of Theorem \ref{thm2}.}
The proof follows  the same line as the proof of Theorem \ref{thm1}.

Assume first that $C^2=L_p(2)<\infty$ and $p/n\leqslant y$ for some $y>0$. Define $X_p^\top AX_p$ and
$X_p^\top BX_p$ in the same way as in \eqref{def}. Then, by Lemma \ref{lm2} (with $\varphi=1/(3b)$),
\[\e_{k-1}\Delta_k\geqslant 1-4C^2\varphi ,\quad 1\leqslant k\leqslant n.\]
 Taking $\varphi=\sqrt{y}/(2 C)$ in \eqref{star}, we get $p/(n\varphi)\leqslant y/\varphi=
2C\sqrt{y}$ and
\[\lambda_p(n^{-1}\bX_{pn}\bX_{pn}^\top)\geqslant 1-4C\sqrt{y}+\frac{CZ}{\sqrt{n}},\]
where
  \[Z=\frac{1}{C\sqrt{n}}\sum_{k=1}^n(\Delta_k-\e_{k-1}\Delta_k).\]
As in the proof of Theorem \ref{thm1}, it follows from Lemma \ref{lm2} that $\e_{k-1}\Delta_k^2\leqslant L_p(2)=C^2$, $1\leqslant k\leqslant n$. Therefore, by Lemma \ref{lm3}, 
$\p(Z<-t)\leqslant \exp\{-t^2/2\},$ $t>0.$

Finally, consider the case with $K_p>0$ ( the case with $K_p=0$ is trivial).  By Lemma \ref{lm2} with $b=(3\varphi)^{-1}$ and $\varphi=1/4$,
\[\e_{k-1}\Delta_k\geqslant \frac{K_p^2}{1+4\varphi}=\frac{K_p^2}{2},\quad 1\leqslant k\leqslant n.\]
 Taking $p/n\leqslant y=K_p^2/16$ in \eqref{star}, we get 
\[\lambda_p(n^{-1}\bX_{pn}\bX_{pn}^\top)\geqslant \frac{K_p^2}{4}+\frac{\sqrt{C_p(4/3)} Z}{\sqrt{n}}\]
for some $Z$ with $\p(Z<-t)\leqslant\exp\{-t^2/2\}$, $t>0$ (see the end of the proof of Theorem \ref{thm1}).
Since $C_p(4/3)\leqslant 4/3$, the variable \[Z_0=\frac{\sqrt{C_p(4/3)}}{\sqrt{ 4/3}} Z\] satisfies  $\p(Z_0<-t)\leqslant\exp\{-t^2/2\}$, $t>0$.  Replacing $Z$ by $Z_0$, we get the result.
\\\\
{\bf Proof of Proposition \ref{p1}.}
If $U$ is non-negative random variable with $\e U=1$, then 
\[\e\min\{U,a\}=\e U-\e(U-a)\bI(U>a)\geqslant 1-\e U\bI(U>a)\geqslant 1-\frac{\e U^{1+\alpha/2}}{a^{\alpha/2}},\]
\[\e\min\{U,a\}=\e U-\int_a^\infty \p(U>t)\,dt\geqslant 1-\int_a^\infty \frac{M}{t^{1+\alpha/2}}\,dt\geqslant 1-\frac{2M}{\alpha a^{\alpha/2}},\]
\[\e U\min\{U,a\}\leqslant
\e U^{1+\alpha/2} a^{1-\alpha/2},\]
\begin{align*}
\e U\min\{U,a\}\leqslant & a\e (U-a)I(U>a)+a^2\p(U>a)+\e \min\{U^2,a^2\}\\ 
&\;\;= a\int_a^\infty \p(U>t)\,dt+a\p(U>a)+\int_0^{a^2}\p(U^2>t)\,dt\\ 
\leqslant&a\int_a^\infty \frac{M}{t^{1+\alpha/2}}\,dt+Ma^{1-\alpha/2}+\int_0^{a^2}f(t,\alpha)\,dt\\
\leqslant& (1+2/\alpha)Ma^{1-\alpha/2}+\begin{cases}
2Ma^{1-\alpha/2}/(1-\alpha/2),&\alpha\in(0,2),\\
2M\log \max\{a,1\}+1,&\alpha=2,\end{cases}
\end{align*}
where $M=\sup\{t^{1+\alpha/2}\p(U>t):t>0\},$ $f(t,\alpha)=Mt^{-1/2-\alpha/4}$ for $\alpha\in(0,2)$ and
\[f(t,2)=\begin{cases}Mt^{-1},&t>1,\\
1,& t\in[0,1].\end{cases}\] 

Putting $U=(X_p,v)^2$ for given $v\in\bR^p$ with $\|v\|=1$ and taking the infimum or the supremum over such $v$ in the above inequalities, we finish the proof.
\\
{\bf Proof of Corollary \ref{cl1}.}
Consider the case $\alpha\in(0,2).$ Set $L=L_p(\alpha)$ and $y=p/n$. By Proposition \ref{p1},
\[c_p(a)-\frac{C_p(a)}{a}\geqslant 1-\frac{2L}{a^{\alpha/2}}\quad \text{and}\quad C_p(2a)\leqslant L\,(2a)^{1-\alpha/2}\leqslant 2La^{1-\alpha/2}.\]
By Theorem \ref{thm1}, 
\begin{align*}
\p(\lambda_p(n^{-1}\X_{pn}\X_{pn}^\top)<1-4La^{-\alpha/2}-5ay)\leqslant& 
\p\big(\sqrt{C_p(2a)}Z/\sqrt{n}<-2La^{-\alpha/2}\big)\\
\leqslant& \p(\sqrt{2La^{1-\alpha/2}}Z/\sqrt{n}<-2La^{-\alpha/2})\\
\leqslant&  \exp\{-L a^{-1-\alpha/2}n\}.
 \end{align*}
 Taking $y=La^{-1-\alpha/2}$, we get the desired inequality.

Consider the case $\alpha=2.$ By Theorem \ref{thm2} with $y=p/n$ and $C=\sqrt{L_p(2)}$, 
\[\p(\lambda_p(n^{-1}\X_{pn}\X_{pn}^\top)<1-(4+\sqrt{2})C\sqrt{y})\leqslant 
\p(CZ/\sqrt{n}<-\sqrt{2}C\sqrt{y})\leqslant \exp\{-yn\}=\exp\{-p\}.\]
\\\\
{\bf Proof of Corollary \ref{cl2}.}  Set $L=L_p(\alpha)$ for given $\alpha\in(0,2)$.  By Proposition \ref{p1},
\[c_p(a)-\frac{C_p(a)}{a}\geqslant 1-\frac{2L}{a^{\alpha/2}}.\]
Therefore, taking in  Theorem \ref{thm1} 
\[a=(4L/\varepsilon)^{2/\alpha}\quad\text{and}\quad p/n\leqslant y=\frac{\varepsilon^{1+2/\alpha}}{10(4L)^{2/\alpha}},\]  we derive the first bound
\[\e \lambda_p(n^{-1}\X_{pn}\X_{pn}^\top)\geqslant 1-\frac{2L}{a^{\alpha/2}}-5ay\geqslant 1-\varepsilon.\]

Similarly, taking $y=\varepsilon^{2}/(16 C^2)$ for $C=\sqrt{L_p(2)}$ in  Theorem \ref{thm2}, we get that
\[\e \lambda_p(n^{-1}\X_{pn}\X_{pn}^\top)\geqslant 1-4C\sqrt{y}\geqslant 1-\varepsilon.\]
\\{\bf Proof of Corollary \ref{cl3}.}
Let $C_0,C_1,C_2>0$ be such that the second bound in Theorem \ref{thm2} holds. Then, for $p/n\leqslant C_2K_p^2,$
\[\p(\lambda_p(n^{ -1}\X_{pn}\X_{pn}^\top)<C_0K_p^2/2)\leqslant 
\p(C_1Z/\sqrt{n}<-C_0K_p^2/2)\leqslant\exp\{-C_0^2K_p^4n/(8C_1^2)\}.\]
Putting $C_0^*=C_0/2,$ $C_1^*=C_0^2/(8C_1^2)$ and $C_2^*=C_2$, we finish the proof.

\section{Appendix}
{\bf Proof of Lemma \ref{l1}.}
By Lemma 2.2 in Srivastava and Vershynin \cite{SV}, if $A-(l+\Delta) I_p\succ 0$ and $q(l+\Delta ,v)/[1+Q(l+\Delta,v)]\geqslant \Delta $, then 
\[\tr(A+vv^\top-(l+\Delta) I_p)^{-1}\leqslant \tr(A-l I_p)^{-1}.\] 
In addition, by Lemma  2.4 in Srivastava and Vershynin \cite{SV}, if $A-lI_p\succ 0$, $\Delta < 1/\varphi$ and $\tr(A-lI_p)^{-1}\leqslant \varphi$, then $A-(l+\Delta )I_p\succ 0$  and
\[\frac{q(l+\Delta ,v)}{1+Q(l+\Delta ,v)}\geqslant \frac{q(l,v)(1-\varphi\Delta)^2}{1+Q(l,v)(1-\varphi\Delta )^{-1}}.\]
Therefore, we only need to show that 
\[\frac{q(l,v)(1-\varphi\Delta)^2}{1+Q(l,v)(1-\varphi\Delta )^{-1}}\geqslant 
\Delta=\frac{q(l,v)}{1+3\varphi q(l,v)+Q(l,v)},\]
since $\Delta\leqslant 1/(3\varphi)$ by construction.

By Bernoulli's inequality, $(1-x)^3\geqslant 1-3x$ whenever $x\in[0,1]$. Hence,
\[\frac{q(l,v)(1-\varphi\Delta)^2}{1+Q(l,v)(1-\varphi\Delta )^{-1}}=
\frac{q(l,v)(1-\varphi\Delta)^3}{1-\varphi\Delta+Q(l,v)}\geqslant 
\frac{q(l,v)(1-\varphi\Delta)^3}{1+Q(l,v)}\geqslant \frac{q(l,v)(1-3\varphi\Delta)}{1+Q(l,v)}=\Delta,\]
where the last equality holds by the definition of $\Delta.$
\\\\{\bf Proof of Lemma \ref{lm1}.}
 We have 
\[\e \frac{U}{1+V}\geqslant \e \frac{\min\{U,a\}}{1+V}\]
for all $a>0.$ By the Cauchy-Schwartz inequality,
\[\e \frac{\min\{U,a\}}{1+V}\e (1+V)\min\{U,a\}\geqslant
 \Big|\e  \frac{\sqrt{\min\{U,a\}}}{\sqrt{1+V}}\,\sqrt{(1+V)\min\{U,a\}}\Big|^2 =|\e\min\{U,a\}|^2.\]
 This gives the first inequality. 
 Tending $a$ to infinity, we get the second inequality.
 
 The last inequality also follows from 
the Cauchy-Schwartz inequality. Namely,
\[\e \frac{U}{1+V} \,\e(1+V)\geqslant \Big|\e \frac{\sqrt{U}}{\sqrt{1+V}}\,\sqrt{1+V}\Big|^2 =|\e\sqrt{U}|^2.\]
\\\\{\bf Proof of Lemma \ref{lm2}.} 
Let $\{v_1,\ldots,v_p\}$ be an orthonormal basis of $\bR^p$ such that \[A=\sum_{i=1}^p a_i v_iv_i^\top\quad\text{and}\quad B=\sum_{i=1}^p b_i v_iv_i^\top,\] where $a_1,\ldots,a_p,b_1,\ldots,b_p>0$ are eigenvalues of $A$ and $B$. Since $\tr A=\sum_{i=1}^pa_i=1$,
$X_p^\top AX_p=\sum_{i=1}^pa_i(X_p,v_i)^2$ and the  function $f(x)=x/(1+c(x+d))$ is concave on $\bR_+$ for any $c,d\geqslant 0$, we have (for $\Delta$ defined in Lemma \ref{lm2})
\[\Delta\geqslant \sum_{i=1}^p a_i\Delta_i\quad\text{for}\quad 
\Delta_i=\frac{(X_p,v_i)^2}{1+b^{-1}((X_p,v_i)^2+X_p^\top B X_p/3)}.\]

Fix $j\in \{1,\ldots,p\}$ and $b>0.$ By Lemma  \ref{lm1}, 
\[\e \Delta_j\geqslant \frac{|\e\min\{(X_p,v_j)^2,a\}|^2}{\e\min\{(X_p,v_j)^2,a\}+b^{-1} C}\quad\text{and}\quad 
\e \Delta_j\geqslant \frac{(\e|(X_p,v_j)|)^2}{1+b^{-1}(1+\tr B/3)}\geqslant \frac{K_p^2}{1+4/(3b)},\]
where $C=\e((X_p,v_j)^2+X_p^\top B X_p/3)\min\{(X_p,v_j)^2,a\}.$ By the second inequality, 
\[\e\Delta\geqslant \sum_{i=1}^pa_i \frac{K_p^2}{1+4/(3b)}=\frac{K_p^2}{1+4/(3b)}.\]

We have $x^2/(x+c)\geqslant x-c$  for all $x,c\geqslant 0$.  This yields that 
\[ \frac{|\e\min\{(X_p,v_j)^2,a\}|^2}{\e\min\{(X_p,v_j)^2,a\}+b^{-1} C}\geqslant
\e\min\{(X_p,v_j)^2,a\}-b^{-1} C.\]

We need to bound $C$ from above.  Obviously, $\e(X_p,v_j)^2\min\{(X_p,v_j)^2,a\}\leqslant C_p(a)$. In addition, since $x\min\{y,a\}\leqslant x\min\{x,a\}+y\min\{y,a\}$ for all $x,y,a\geqslant 0,$  we have 
 \[\e(X_p^\top B X_p)\min\{(X_p,v_j)^2,a\}=
\sum_{i=1}^pb_i\e(X_p,v_i)^2\min\{(X_p,v_j)^2,a\}\leqslant 2\tr B\cdot C_p(a)\leqslant 2 C_p(a).\]
Hence, $C\leqslant  5C_p(a)/3.$ 
Combining all estimates together yields 
\[\e\Delta\geqslant  c_p(a)-\frac{5 C_p(a)}{3b}. \]

Let us now prove that $\e \Delta^2 \leqslant C_p(b)$. We have 
\[\Delta^2\leqslant \frac{(X_p^\top A X_p)^2 }{(1+b^{-1}X_p^\top A X_p)^2}\leqslant
 \frac{(X_p^\top A X_p)^2 }{1+b^{-1}X_p^\top A X_p}.\]
Consider the function $f(x)=x^2/(1+b^{-1}x)$, $x\geqslant 0$. Its derivative
\[f'(x)=\frac{2x}{1+b^{-1}x}-\frac{b^{-1}x^2}{(1+b^{-1}x)^2}=
\frac{2x+b^{-1}x^2}{(1+b^{-1}x)^2}=b \frac{2bx+x^2}{(b+x)^2}=b\Big(1-\frac{b^2}{(b+x)^2}\Big)\]
is increasing on $\bR_+$. This means that $f=f(x)$ is convex and
\[\e \frac{(X_p^\top A X_p)^2 }{1+a^{-1}X_p^\top A X_p}\leqslant \sum_{i=1}^p a_i\e\frac{(X_p,v_i)^4}{1+b^{-1}(X_p,v_i)^2}\leqslant
\sum_{i=1}^p a_i\e(X_p,v_i)^2\min\{(X_p,v_i)^2,b\}.\]
The latter gives the desired inequality $\e\Delta^2\leqslant \tr A\cdot C_p(b)=C_p(b).$

Now consider the case with $L_p(2)<\infty$. By Lemma \ref{lm1},
\[\e \Delta \geqslant  1/[1+b^{-1}(\e (X_p^\top AX_p)^2+ \e (X_p^\top AX_p)(X_p^\top BX_p)/3)].\]
Since the function $f(x)=x^2$ is convex on $\bR$, $X_p^\top AX_p=\sum_{i=1}^n a_i (X_p,v_i)^2$ and $\tr A=1$, we get that
\[ \e (X_p^\top AX_p)^2\leqslant
\sum_{i=1}^n a_i \e (X_p,v_i)^4\leqslant L_p(2).\]
Similarly,
\[ \e (X_p^\top BX_p)^2\leqslant (\tr B)^2 \e\Big(\frac{X_p^\top BX_p}{\tr B}\Big)^2\leqslant L_p(2),\]
where we have used that $\tr B\leqslant 1$. Applying the Cauchy-Schwartz inequality yields that 
\[\e (X_p^\top AX_p)(X_p^\top BX_p)\leqslant\sqrt{\e(X_p^\top AX_p)^2\e (X_p^\top BX_p)^2}\leqslant L_p(2).\] 
To finish the proof, we  only need to note that
\[ 1/[1+b^{-1}(\e (X_p^\top AX_p)^2+ \e (X_p^\top AX_p)(X_p^\top BX_p)/3)]\geqslant 
 \frac{1}{1+4L_p(2)b^{-1}/3}\geqslant 1-\frac{4L_p(2)}{3b}.\]
\\\\{\bf Proof of Lemma \ref{lm3}.} 
Since $e^{-x}\leqslant 1-x+x^2/2$ for all $x\geqslant 0,$ we have 
\begin{align*}
\e(e^{-\lambda D_k}|\cF_{k-1})\leqslant&
1-\lambda\e(D_k|\cF_{k-1})+\frac{\lambda^2\e(D_k^2|\cF_{k-1})}{2}\\
\leqslant &
1-\lambda\e(D_k|\cF_{k-1})+\frac{\lambda^2}{2}\\
\leqslant &\exp\{-\lambda\e(D_k|\cF_{k-1})+\lambda^2/2\}
\end{align*}
for any $\lambda>0.$ Therefore, $\e(e^{-\lambda (D_k-\e(D_k|\cF_{k-1}))}|\cF_{k-1})\leqslant \exp\{\lambda^2/2\}$ and
\begin{align*}
\p\Big(\sum_{k=1}^n (D_k-\e(D_k|\cF_{k-1}))<-t\sqrt{n}\Big)\leqslant & e^{-\lambda t\sqrt{n}}\e\exp\Big\{-\lambda \sum_{k=1}^n(D_k-\e(D_k|\cF_{k-1}))\Big\}\\
\leqslant & \exp\{n\lambda^2/2-\lambda t\sqrt{n}\},
\end{align*}
where the last bound could be obtained iteratively by the law of iterated mathematical expectations.
Putting $\lambda= t/\sqrt{n}$, we derive that $\p(Z<-t)\leqslant \exp\{-t^2/2\}$, $t>0$.




\end{document}